\documentclass[11pt]{amsart}

\usepackage{amssymb,amsmath,enumerate,epsfig}%,epsf}%t1enc
\usepackage{amsfonts,color}
\usepackage{a4,latexsym}
\usepackage{parskip}
\usepackage[all]{xy}
\usepackage{hyperref}
\hypersetup{pdftitle=Enriques}

\newcommand{\C} {\mathbb{C}}
\newcommand{\Q} {\mathbb{Q}}

\newcommand{\F}{\mathbb{F}}
\newcommand{\Z}{\mathbb{Z}}

\newcommand{\PP}{\mathbb{P}}
\newcommand{\NS}{\mathop{\rm NS}}

\newcommand{\MWL}{\mathop{\rm MWL}}
\newcommand{\Hom}{\mathop{\rm Hom}}
\newcommand{\Km}{\mathop{\rm Km}}

\newtheorem{Theorem}{Theorem}[section]
\newtheorem{Proposition}[Theorem]{Proposition}
\newtheorem{Lemma}[Theorem]{Lemma}

\theoremstyle{remark}
\newtheorem{Remark}[Theorem]{Remark}

\newtheorem{Example}[Theorem]{Example}

\theoremstyle{definition}

\begin{document}

\title{A note on the supersingular K3 surface
of Artin invariant 1}

\author{Matthias Sch\"utt}
\address{Institut f\"ur Algebraische Geometrie, Leibniz Universit\"at
  Hannover, Welfengarten 1, 30167 Hannover, Germany}
\email{schuett@math.uni-hannover.de}
\urladdr{http://www.iag.uni-hannover.de/$\sim$schuett/}

%\subjclass[2000]{14J28; 14J27, 14J50, 14F22}
%
%\keywords{Enriques surface, K3 surface, elliptic fibration, Mordell-Weil group, automorphism, Brauer group}
%
%
%\thanks{Partial funding from DFG grant Hu 337/6-1 is gratefully acknowledged}
%
%\dedicatory{Dedicated to Tetsuji Shioda on the occasion of his 70th birthday}
%

\date{\today}

 \begin{abstract}
We prove that  the supersingular K3 surface
of Artin invariant 1 in characteristic $p$ (where $p$ denotes an arbitrary prime)
 admits a model over $\F_p$ with Picard number 21.
 \end{abstract}
 
 \maketitle

\section{Introduction}

This note concerns an explicit problem about the supersingular K3 surface $X$
of Artin invariant 1 in characteristic $p$
where $p$ denotes an arbitrary prime.
Over $\bar\F_p$ this surface is unique
by \cite{Ogus}
(see \cite{RS2} for characteristic $2$).
Here we prove the existence of a  model over $\F_p$ with ideal properties:

\begin{Theorem}
\label{thm}
The K3 surface $X$ admits a model over $\F_p$ with Picard number 21.
\end{Theorem}

The Picard number is the rank of the N\'eron-Severi group $\NS(X)$
consisting of divisors up to algebraic equivalence 
(or, in the context of K3 surfaces, linear or numerical equivalence).
The above model attains the geometric Picard number 22 over the quadratic extension $\F_{p^2}$.
We remark that this maximum cannot be attained over $\F_p$, 
or in fact over any finite field $\F_{p^e}$ with $e$ odd,
by \cite[(6.8)]{Artin} 
(see also \cite[Thm.~4.4]{S-MJM}).
In this sense, Theorem \ref{thm} represents the ideal situation over $\F_p$.

This note is organised as follows.
The next section gives a motivation for the problem
-- in fact a solution for most characteristics.
Then Section \ref{s:proof} gives a general proof
that builds on Shioda--Inose structures
and a sandwich picture developed by Shioda in \cite{Sandwich}, \cite{Sh-Murre}
and extended in \cite{HSa}.

\section{Motivation: singular K3 surfaces}

The problem of Theorem \ref{thm} arose from discussions with H.~Ohashi
who considered a very specific elliptic fibration 
on the supersingular K3 surface $X$
of Artin invariant 1 in characteristic $11$ in \cite{Ohashi}.
The author's initial idea was that the fact 
that $X$ has a model with $\NS$ fully defined over $\F_{11^2}$
might simplify some arguments in \cite{Ohashi}.
Subsequently this led to the general statement of Theorem \ref{thm}.

A basic approach to see the claim in characteristic $11$ is reduction from characteristic zero
(see Example \ref{ex}).
Indeed this provides a convenient way to produce supersingular K3 surfaces.
Here we use as input singular K3 surfaces,
i.e.~complex K3 surfaces whose Picard number attains Lefschetz' bound of $h^{1,1}=20$.
These can be classified completely in terms of their transcendental lattices by \cite{SI}.
In particular, each one is defined over some number field.
However, over all number fields of degree not exceeding some given bound,
there are only finitely many singular K3 surfaces up to $\bar\Q$-isomorphism (see \cite{S-NS}).
Note the following subtlety which is in contrast to the result of Theorem \ref{thm}:
there are singular K3 surfaces over arbitrarily large number fields
whose moduli point is defined over $\Q$.
By this rough statement we mean that the surface does not admit a model over any smaller field,
yet it is $\bar\Q$-isomorphic to all its Galois-conjugates.
In comparison, Theorem \ref{thm} states that the $\bar\F_p$-moduli point formed by the supersingular K3 surface $X$ corresponds indeed to a K3 surface over the prime field $\F_p$.

Specifically consider singular K3 surfaces over $\Q$ with all of the N\'eron-Severi group defined over $\Q$ as well.
By \cite{S-NS} there are 13 such surfaces up to $\bar\Q$-isomorphism.
These are in 1-to-1 correspondence with elliptic curves over $\Q$ with complex multiplication.
Such a singular K3 surface $X$ gives rise to the discriminant $d<0$ of the N\'eron-Severi group,
and thus to an imaginary quadratic field $K=\Q(\sqrt{d})$.

\begin{Lemma}
\label{fact}
The reduction of $X$ modulo some good prime $p\nmid d$ is supersingular
if and only if $p$ is inert in $K/\Q$.
\end{Lemma}

\begin{Example}
\label{ex}
The supersingular K3 surface with Artin invariant 1 in characteristic $11$
(as studied in \cite{Ohashi})
arises from the singular K3 surfaces with discriminant $-3$ or $-4$ by reduction.
\end{Example}

Lemma \ref{fact} can be seen geometrically by way of Shioda-Inose structures
as we will exploit in Section \ref{s:proof}.
Alternatively one can argue with modularity after Livn\'e \cite{L}.
Since the corresponding Hecke eigenform has weight 3 and nebentypus character $\chi$ associated to $K$,
one finds the characteristic polynomial of Frobenius on $H_\text{\'{e}t}^2(X\otimes\bar\F_p, \Q_\ell)$
as
\[
P(X,T) = (T-p)^{20} (T^2-a_pT+\chi(p)p^2) =  (T-p)^{21} (T+p).
\]
Here first equality holds generally with $a_p$ denoting the eigenvalues of the eigenform
and the 20-fold factor $(T-p)$ coming from $\NS$ in characteristic zero.
Meanwhile the second equality depends on the choice of an inert prime.
Since the Tate conjecture is known for elliptic K3 surfaces with section by \cite{ASD}, 
and every singular K3 surface has such a fibration (as exploited in Section \ref{s:proof})
we infer that the reduction of $X$ is supersingular; 
more precisely we find that $X$ has Picard number 21 over $\F_p$ and 22 over $\F_{p^2}$.
It remains to compute the Artin invariant of the supersingular reduction $X$.
For this we refer to a result of Shimada \cite[Proposition 1.0.1]{Shimada}:

\begin{Proposition}
The reduction of a singular K3 surface modulo a supersingular prime
has Artin invariant 1.
\end{Proposition}

In conclusion, the supersingular K3 surfaces $X$ derived as above do exactly fit
with Theorem \ref{thm}.
It is only the primes which split
% or ramify, but that does not occur
 in each imaginary quadratic field of class number one
where the above construction fails to produce the required models.
In order to cover these primes ($2^{-9}$th of all primes, the first one being $15073$)
we will employ a geometric argument along the lines of \cite{HSa} in the next section.

\section{Geometric approach}
\label{s:proof}

Another prototype of (supersingular) K3 surfaces are Kummer surfaces.
Here we start with an abelian variety $A$,
quotient by the involution and desingularise to obtain a K3 surface $\Km$.
If our initial abelian variety is a product of two elliptic curves $E, E'$,
then $\Km$ is supersingular if and only if both $E$ and $E'$ are.
However, these Kummer surfaces are not sufficient for our purpose 
of proving Theorem \ref{thm}.
Namely they inherit too big a Galois action on $\NS$ from the abelian variety.
Possibly this can also stem from the 2-torsion points of $E$ and $E'$,
but always from the cohomology of $A$ since the characteristic polynomial 
of Frobenius on $\wedge^2 H_\text{\'et}^1(A\otimes\bar\F_p,\Q_\ell)\subset H_\text{\'et}^2(A\otimes\bar\F_p,\Q_\ell)$ implies $\rho(E\times E'/\F_p)=4$ so that
the natural model of the Kummer surface has $\rho(\Km(E\times E')/\F_p)\leq 20$.

It is instructive to note the parallel
that complex Kummer surfaces do also not suffice to treat all singular K3 surfaces.
Continuing the analogy, we will invoke the concept of Shioda--Inose structures
in order to give a complete proof of Theorem \ref{thm}.

\subsection{Shioda--Inose and sandwich structure}

Over $\C$, the notion of Shioda--Inose structure refers to a pair of an abelian surface $A$ and a K3 surface $X$ 
(not necessarily singular)
with the same transcendental lattice
such that $X$ admits a rational map of degree two to $\Km(A)$:
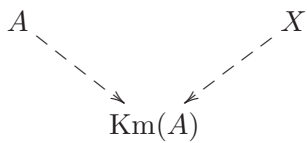
 \[
  \xymatrix{A \ar@{-->}[dr] && X\ar@{-->}[dl]\\
 & \Km(A)&}
 \]
In \cite{SI}, Shioda--Inose prove that any singular K3 surface $X$ fits into such a structure
and can thus be described in terms of products of isogenous CM-elliptic curves.
In \cite{Sandwich} Shioda extends this construction to show
that $X$ is in fact sandwiched by $\Km(A)$.
In particular this implies that $X$ and $\Km(A)$ have the same geometric Picard number
regardless of the characteristic
(a fact that would follow over $\C$ from a notion of isogeny due to Inose \cite{Inose}).

All these constructions are exhibited in terms of explicit algebraic equations (which we give below),
so they directly apply to characteristic $p$.
Here we only have to take extra care of characteristics $2$ and $3$ 
where the elliptic fibrations involved degenerate.
However, for those characteristics, explicit models guaranteeing Theorem \ref{thm}
have been exhibited, for instance, in \cite{S-MJM} and \cite{ST},
so we shall omit them in the sequel.

The next diagram gives a brief schematic sketch of some of the elliptic fibrations in the product case:
$$
\begin{array}{ccccccccc}
&&&&&& \Km(E\times E') &&\\
&&&&& \swarrow && \searrow \pi_1 &\\
A &&&& X &&&& \PP^1\\
\downarrow & \searrow && \swarrow && \searrow\pi_X && \swarrow &\\
E & & \Km(E\times E') &&&& \PP^1 &&\\
& \searrow & \phantom{\pi}\downarrow\pi_0 &&&&&&\\
&& \PP^1 &&&&&&
\end{array}
$$
We continue by explaining how these fibrations arise 
over some field $k$ of characteristic different from $2$.
Throughout we work with Weierstrass models 
\begin{eqnarray}
\label{eq:WF-E}
 E:\;\; y^2 = f(x),\;\;\;\; E': \;\; y^2 = g(x)
\end{eqnarray}
with cubic polynomials $f, g\in k[x]$.
Then $\Km(E\times E')$ admits a birational model
\begin{eqnarray}
\label{eq:Km}
 \Km(E\times E'):\;\; f(t) y^2 = g(x).
\end {eqnarray}
In terms of this model the elliptic fibrations are as follows:
\begin{enumerate}
\item[1.]
The elliptic fibration $\pi_0$ is given by projection onto the projective $t$-line.
This fibration is always isotrivial with four singular fibres of type $I_0^*$ in Kodaira's notation.
%
%with the structure of an elliptic curve over the function field $H(d)(t)$.
%We denote the corresponding elliptic fibration by the pair $(X', \pi)$.
%This fibration has singular fibres of type $I_0^*$ at $\infty$ and at the zeroes of $f(t)$.
%Over $\bar\Q$ we have $\MW(X',\pi)=\Z^2\times (\Z/2\Z)^2$ with torsion sections given by the roots of $g(x)$.
\item[2.]
The elliptic fibration $\pi_1$ is given in terms of \eqref{eq:Km} by projection onto the $y$-line.
If $E\not\cong E'$, then this fibration has only two reducible fibres, both of Kodaira type $IV^*$.
\item[3.]
Finally $X$ arises from $\Km(E\times E')$ as the quotient by the involution
which composes $y\mapsto -y$ with the hyperelliptic involution of fibration $\pi_1$
(once a zero section is chosen).
Clearly this induces an elliptic fibration $\pi_X$ on $X$
which has two type $II^*$ fibres if $E\not\cong E'$.
(This is often referred to as Inose's fibration.)
\end{enumerate}
The next step consists in specialising to the situation where $E\cong E'$. 
Unless $j(E)=0, 12^3$, 
the fibration $\pi_X$ attains exactly one additional reducible fibre which has type $I_2$.
This fibre is duplicated on the fibration $\pi_1$.

Since the two special j-invariants $j(E)=0,12^3$ can be understood completely 
by reducing the corresponding singular K3 surfaces over $\Q$ (of discriminant $-3, -4$),
we shall exclude these cases in the sequel without further mention.

\subsection{Proof of Theorem \ref{thm}}
\label{ss:proof}

From now on, we fix the prime $p$.
The overall idea is to pick some supersingular elliptic curve $E$ over $\F_p$
and consider the K3 surface $X$ arsing from the Shioda--Inose structure for $E\times E$.
Such a curve is characterised by its trace being zero, or equivalently
$\# E(\F_p)=p+1$.
Thus its existence over $\F_p$ follows, for instance, from Honda's theorem \cite{Honda}
which states that there exists an elliptic curve over $\F_p$ with any trace $a\in\Z,
|a |\leq 2\sqrt{p}$.

The elliptic fibration $\pi_X$ thus obtained has Mordell-Weil rank 3 over $\bar\F_p$.
It seems feasible to apply abstract lattice theoretic arguments along the lines of \cite[\S 3]{HSa}, 
combined with lifting to a singular K3 surface, to analyse the possible Galois action on $\MWL(X,\pi_X)$ 
and produce a quadratic twist, if necessary, with MW-rank 2 over $\F_p$.
However, we decided to pursue a direct geometric approach. % that seems more direct.

A detailed analysis of the Mordell-Weil groups of the elliptic fibrations in question
has been carried out over algebraically closed fields by Shioda in \cite{Sh-Murre}.
As in \cite[\S 4]{HSa}, we throw in just a little bit of extra thought to make the argument work
over a non-algebraically closed field $k$.
%Throughout we consider some finite field $k$ of characteristic $p$ which we suppress in the notation. 
Later on we will specialise to the case $k=\F_p$.

Consider the lattice $\Hom(E, E)$ endowed with a norm given by the degree.
By  \cite[Prop.~3.1]{Sh-Murre}, there is an isomorphism of lattices
\begin{eqnarray}
\label{eq:Hom}
 \Hom(E, E) \cong \MWL(\Km(E\times E),\pi_0).
\end{eqnarray}
By construction this is clearly Galois-equivariant.
Hence, for $E/\F_p$ supersingular, we deduce from the argument given at the beginning of this section,
that $\pi$ has MW-rank 2 over $\F_p$ and 4 over $\F_{p^2}$.

Our aim is to compare these lattices to a certain sublattice of $\NS(\Km(E\times E))$
arising from fibration $\pi_1$.
For this purpose consider $\pi_X$ for the moment.
The singular fibres of type $II^*$ together with the zero section
generate the unimodular sublattice $U+E_8(-1)^2\subset\NS(X)$
which is completely defined over $k$.
We denote its orthogonal complement by $L$.
Usually the lattice $L$ is exactly the Mordell-Weil lattice of $\pi_X$ up to sign,
but if $E\cong E'$, as is presently the case, $L(-1)$ is comprised of the root lattice $A_1$
corresponding to an additional reducible fibre, 
and the Mordell-Weil lattice of rank $3$ over $\bar k$
(assuming that $E$ is supersingular). %, but $j(E)\neq 0, 12^3$).

We can pull-back the lattice $L$ via the quotient map $\Km(E\times E)\dasharrow X$
to obtain a natural sublattice $L(2)\subset\NS(\Km(E\times E))$.
The crucial point in Shioda's argument in \cite[\S 4  \& 7]{Sh-Murre}
is a geometric transition $\bar k$  between the fibrations $\pi_0$ and $\pi_1$ on $\Km(E\times E)$
which implies that 
\[
L(2) \cong \Hom(E, E)(4).
\]
In order to guarantee the Galois-equivariance of this lattice isomorphism,
it suffices that the fibration $\pi_1$ given by \eqref{eq:Km} 
has a base point over $k(y)$, i.e.~$E$ has a 2-torsion point giving a zero of $f=g$.
Recall that presently we are concerned with a supersingular elliptic curve $E$ over $\F_p$.
Any such curve has an $\F_p$-rational 2-torsion point  (outside characteristic 2)
 simply because, the trace being zero, $\# E(\F_p)=p+1$ is even.
 
% Assumption throughout: $j(E)\neq 0, 12^3$ -- else reduction from singular K3 with $\NS$ fully over $\Q$
 
 In conclusion we find for the given model of $X$ over $\F_p$
 that $L$ has rank 2 over $\F_p$ and rank 4 over $\F_{p^2}$.
 Since both components of the $I_2$ fiber are clearly defined over $\F_p$,
 this implies Mordell-Weil rank 1 over $\F_p$ and rank 3 over $\F_{p^2}$.
 But then the quadratic twist $X'$ of $X$ over $\F_{p^2}$ automatically has reversed MW-ranks:
 2 over $\F_p$ and 3 over $\F_{p^2}$.
Summing up, $X'$ gives an alternative model of $X$ defined over $\F_p$
with $\rho(X'/\F_p)=21$. 
This verifies Theorem \ref{thm}. \qed

\begin{Remark}
Implicitly the above argument uses the fact
that quadratic twists do not affect fibers of type $II^*$ and $I_2$
(such as on fibration $\pi_X$ on $X$).
This does not hold for fibers of type $IV^*$, for instance,
so quadratic twisting has a fundamentally different effect on fibration $\pi_1$ on $\Km(E\times E)$.
\end{Remark}

\subsection*{Acknowledgements}
This note grew out of a discussion with H.~Ohashi.
The author thanks T.~Shioda and S.~Wewers  for helpful comments.

\end{document}